\documentclass[12pt]{article}
\usepackage{amsmath}
\usepackage{amssymb}
\usepackage{vatola}

\def\q{\quad}
\def\qq{\qquad}
\def\qtq#1{\q\t{#1}\q}
\def\mod#1{\ (\text{\rm mod}\ #1)}
\def\t{\text}
\def\f{\frac}
\def\e{\equiv}
\def\b{\binom}
\def\sls#1#2{(\f{#1}{#2})}
 \def\ls#1#2{\big(\f{#1}{#2}\big)}
\def\Ls#1#2{\Big(\f{#1}{#2}\Big)}
\let \pro=\proclaim
\let \endpro=\endproclaim

\begin{document}
\leftline{Publ. Math. Debrecen 96(2020), no.3-4, 315-346.}
\par\q
 \centerline {\bf
Congruences involving binomial coefficients and  Ap\'ery-like
numbers}
\par\q\newline
\centerline{Zhi-Hong Sun}
 \abstract{For $n=0,1,2,\ldots$ let $W_n=\sum_{k=0}^{[n/3]}\binom{2k}k
 \binom{3k}k\binom n{3k}(-3)^{n-3k}$,
 where $[x]$ is the greatest integer not exceeding $x$. Then
 $\{W_n\}$ is an Ap\'ery-like  sequence. In this paper we deduce many congruences
 involving $\{W_n\}$, in particular we determine
$\sum_{k=0}^{p-1}\binom{2k}k\frac{W_k}{m^k}\pmod p$ for
$m=-640332,-5292,-972,-108,-44,-27,-12,8,54,243$ by using binary
quadratic forms, where $p>3$ is a prime. We also prove several
congruences for generalized Ap\'ery-like numbers, and pose 29
challenging conjectures on congruences involving binomial
coefficients and Ap\'ery-like numbers.
 \par\q
 }
 \endabstract
\let\thefootnote\relax \footnotetext {The author is supported by
the National Natural Science Foundation of China (Grant No.
11771173).}
\section*{1. Introduction}
\par For $s>1$ let $\zeta(s)=\sum_{n=1}^{\infty}\f 1{n^s}$.
In 1979, in order to prove $\zeta(2)$ and
 $\zeta(3)$ are irrational,
 Ap\'ery [4] introduced the Ap\'ery numbers $\{A_n\}$ and $\{A'_n\}$
 given by
 $$A_n= \sum_{k=0}^n\binom nk^2\binom{n+k}k^2\qtq{and}A'_n=\sum_{k=0}^n\b
 nk^2\b{n+k}k.$$
It is well known (see [7]) that
$$\align &(n+1)^3A_{n+1}=(2n+1)(17n(n+1)+5)A_n-n^3A_{n-1}\q (n\ge 1),
\\&(n+1)^2A'_{n+1}=(11n(n+1)+3)A'_n+n^2A'_{n-1}\q (n\ge 1).\endalign$$
\par
Let $\Bbb Z$ and $\Bbb Z^+$ be the set of integers and the set of
positive integers, respectively.  The first kind of Ap\'ery-like
numbers $\{u_n\}$ satisfies
$$u_0=1, \ u_1=\ b,\ (n+1)^3u_{n+1}=(2n+1)(an(n+1)+b)u_n-cn^3u_{n-1}\q (n\ge 1),\tag 1.1$$
where $a,b,c\in\Bbb Z$ and $c\not=0$. Let $[x]$ be the greatest
integer not exceeding $x$, and let
$$\align
&D_n=\sum_{k=0}^n\b nk^2\b{2k}k\b{2n-2k}{n-k},\q T_n=\sum_{k=0}^n\b
nk^2\b{2k}n^2, \\& b_n=\sum_{k=0}^{[n/3]}\b{2k}k\b{3k}k\b
n{3k}\b{n+k}k(-3)^{n-3k}.
\endalign$$
Then $\{A_n\}$, $\{D_n\}$, $\{b_n\}$ and $\{T_n\}$  are the first
kind of Ap\'ery-like numbers with
$(a,b,c)=(17,5,1),(10,4,64),(-7,-3,81),(12,4,16)$, respectively. The
numbers $\{D_n\}$ are called Domb numbers, and $\{b_n\}$ are called
Almkvist-Zudilin numbers. For $\{A_n\}$, $\{D_n\}$, $\{b_n\}$ and
$\{T_n\}$ see A005259, A002895, A125143 and A290575 in Sloane's
database ``The On-Line Encyclopedia of Integer Sequences" [40]. For
the congruences involving $T_n$ see the author's recent paper [33].
\par In 2009, Zagier [41] studied the second kind of Ap\'ery-like numbers
$\{u_n\}$ given by
$$u_0=1,\ u_1=b\qtq{and}(n+1)^2u_{n+1}=(an(n+1)+b)u_n-cn^2u_{n-1}\ (n\ge 1),\tag 1.2$$
where $a,b,c\in\Bbb Z$ and $c\not=0$.  Let
 $$\align &f_n=\sum_{k=0}^n\b nk^3=\sum_{k=0}^n\b nk^2\b{2k}n,\q
 \\&S_n=\sum_{k=0}^{[n/2]}\b{2k}k^2\b n{2k}4^{n-2k}=\sum_{k=0}^n\b nk\b{2k}k\b{2n-2k}{n-k},
 \\&a_n=\sum_{k=0}^n\b nk^2\b{2k}k,
 \q Q_n=\sum_{k=0}^n\b nk(-8)^{n-k}f_k,
 \\&W_n=\sum_{k=0}^{[n/3]}\b{2k}k\b{3k}k\b n{3k}(-3)^{n-3k}.\endalign$$
  In [41] Zagier stated that $\{A'_n\},\ \{f_n\},\ \{S_n\},\ \{a_n\},\
 \{Q_n\}$
and $\{W_n\}$ are the second kind of Ap\'ery-like sequences with
$(a,b,c)=(11,3,-1),(7,2,-8),(12,4,32),(10,3,9),(-17,-6,72),$ $
(-9,-3,27)$, respectively. The  sequence
 $\{f_n\}$ is called Franel numbers.
 In [14,31,32]  the author
 systematically investigated identities and congruences for sums
 involving $S_n$ or $f_n$.
  For $\{A'_n\},\ \{f_n\},\ \{S_n\},\ \{a_n\},\
 \{Q_n\}$
and $\{W_n\}$, see A005258, A000172, A081085, A002893, A093388 and
A291898 in Sloane's database [40].
\par Ap\'ery-like numbers have fascinating properties and they are
concerned with Picard-Fuchs differential equation, modular forms,
hypergeometric series, elliptic curves, series for $\f 1{\pi}$,
supercongruences, binary quadratic forms, combinatorial identities,
Bernoulli numbers and Euler numbers. See typical papers
[1,2,3,6,8,9,10,11,13,17,18,21,30,35,36,39].
\par  For $a\in\Bbb Z$ and given odd prime $p$ let $\sls ap$ be the
Legendre symbol. For a prime $p$ let $\Bbb Z_p$ be the set of
rational numbers whose denominator is not divisible by $p$. For
positive integers $a,b$ and $n$, if $n=ax^2+by^2$ for some integers
$x$ and $y$, we briefly write that $n=ax^2+by^2$.

\par In Section 2 we obtain some congruences for sums involving $W_n$.
 From [9, (6.4)] we know
that
$$\Big(\sum_{k=0}^{\infty}W_kx^k\Big)^2=\f 1{1-27x^2}
\sum_{k=0}^{\infty}\b{2k}k\Big(\f{x(1+9x+27x^2)}{(1-27x^2)^2}\Big)^kW_k.\tag
1.3$$ We prove the $p-$analogue of (1.3):
$$\Big(\sum_{k=0}^{p-1}W_kx^k\Big)^2
\e\sum_{k=0}^{p-1}\b{2k}k\Ls{x(1+9x+27x^2)}{(1-27x^2)^2}^kW_k\mod
p,\tag 1.4$$ where $p>3$ is a prime, $x\in\Bbb Z_p$ and
$(x+3)(1+9x+27x^2)(1+9x)(1+27x^2)(1-27x^2)\not\e 0\mod p$.
\par
Suppose that $p>3$ is a prime, $n\in\Bbb Z_p$ and $n(n-12)\not\e
0\mod p$. In Section 2, we show that
$$\sum_{k=0}^{p-1}\b{2k}k\f{W_k}{(n-12)^k}
\e\Ls {n(n-12)}p\sum_{k=0}^{p-1}\f{\b{2k}k\b{3k}k\b{6k}{3k}}
{n^{3k}}\mod p.\tag 1.5$$ As consequences, we determine
$\sum_{k=0}^{p-1}\b{2k}k\f{W_k}{m^k}\mod p$ for
$m=-640332,-5292,-972,-108,$ $-44,-27,8,54,243$ by using binary
quadratic forms. We also determine $W_{\f{p-1}2}$,
$\sum_{k=0}^{p-1}\f{W_k}{(-3)^k}$, $\sum_{k=0}^{p-1}\f{W_k}{(-9)^k}$
and $\sum_{k=0}^{p-1}\b{2k}k\f{W_k}{(-12)^k}$ modulo $p$.

\par Let $p$ be an odd prime. In 2010 Jarvis and Verrill [13] discovered relations between $u_n$
and $u_{p-1-n}$ modulo $p$ for $u_n=A_n',a_n,f_n$ or $S_n$. For
example, they proved $f_n\e (-8)^nf_{p-1-n}\mod p$ for
$n=0,1,\ldots,p-1$. In Section 3 we establish a vast generalization
of such congruences
 for generalized Ap\'ery-like numbers $\{u_n\}$. See
Theorem 3.1.
\par  In [20] S. Ramanujan made some conjectures for $1/\pi$, which
involve the following four sums
$$\sum_{k=0}^{\infty}(ak+b)\f{\b{2k}k^3}{m^k},\
\sum_{k=0}^{\infty}(ak+b)\f{\b{2k}k^2\b{3k}k}{m^k},\
\sum_{k=0}^{\infty}(ak+b)\f{\b{2k}k^2\b{4k}{2k}}{m^k},\
\sum_{k=0}^{\infty}(ak+b)\f{\b{2k}k\b{3k}k\b{6k}{3k}}{m^k}.$$
 The author's brother Z.W. Sun in [39] and the author in [25] posed many
conjectures on congruences for
$$\sum_{k=0}^{p-1}\f{\b{2k}k^3}{m^k},\q
\sum_{k=0}^{p-1}\f{\b{2k}k^2\b{3k}k}{m^k},\q
\sum_{k=0}^{p-1}\f{\b{2k}k^2\b{4k}{2k}}{m^k},\q
\sum_{k=0}^{p-1}\f{\b{2k}k\b{3k}k\b{6k}{3k}}{m^k}\mod{p^2},\tag
1.6$$
 where $m\in\Bbb Z$ and $p$ is an odd prime with $p\nmid m$.
Some of such conjectures were proved by the author in [27-29]. In
particular, most of conjectures were solved when the modulus is $p$.
Let $p>3$ be a prime. Recently Liu [16] conjectured congruences for
$$\sum_{k=0}^{p-1}\f{\b{2k}k^2\b{3k}k}{108^k},\q
\sum_{k=0}^{p-1}\f{\b{2k}k^2\b{4k}{2k}}{256^k},\q
\sum_{k=0}^{p-1}\f{\b{2k}k\b{3k}k\b{6k}{3k}}{1728^k}\mod {p^3}$$ in
terms of $p-$adic gamma functions.  Based on calculations with
Maple, in Section 4 we pose 29 challenging conjectures on
congruences involving sums in (1.6) or Ap\'ery-like numbers. See
Conjectures 4.1-4.29.

\section*{2. Congruences for sums involving $W_n$}
\par For any nonnegative integer $n$, define
$$W_n(x)=\sum_{k=0}^{[n/3]}\b{2k}k\b{3k}k\b n{3k}x^{n-3k}.$$
Then $W_n=W_n(-3).$ In this section we establish some congruences
involving $W_n$ and $W_n(x)$ modulo a prime.
 We begin with three useful lemmas.
 \pro{Lemma 2.1} Let $n$ be a nonnegative
integer. Then
$$\sum_{k=0}^n\b nkW_k(x)y^{n-k}=W_n(x+y).$$
\endpro
Proof. It is clear that
$$\align \sum_{k=0}^n\b nkW_k(x)y^{n-k}&=\sum_{k=0}^n
\b nky^{n-k}\sum_{r=0}^k\b{2r}r\b{3r}r\b k{3r}x^{k-3r}
\\&=\sum_{r=0}^n\b{2r}r\b{3r}ry^{n-3r}
\sum_{k=r}^n\b nk\b k{3r}\Ls xy^{k-3r}
\\&=\sum_{r=0}^n\b{2r}r\b{3r}ry^{n-3r}\sum_{k=3r}^n\b
n{3r}\b{n-3r}{k-3r}\Ls xy^{k-3r}
\\&=\sum_{r=0}^n\b{2r}r\b{3r}r\b n{3r}y^{n-3r}
\Big(1+\f xy\Big)^{n-3r}=W_n(x+y).
\endalign$$ This proves the lemma.
\par\q
\par
Let $\{P_n(x)\}$ be the famous Legendre polynomials given by
$$\align P_n(x)&=\sum_{k=0}^n\b nk\b{n+k}k\Ls{x-1}2^k
=\f 1{2^n}\sum_{k=0}^{[n/2]}\b
nk(-1)^k\b{2n-2k}nx^{n-2k}.\endalign$$
 It is well known that
 $$P_0(x)=1,\ P_1(x)=x,\ (n+1)P_{n+1}(x)=(2n+1)xP_n(x)-nP_{n-1}(x)\q(n\ge 1).$$
\pro{Lemma 2.2} Suppose that $p>3$ is a prime and $m,x\in\Bbb Z_p$
with $mx\not\e 0\mod p$. Then
$$\align \sum_{k=0}^{p-1}\f{W_k(x+m)}{m^k}&\e W_{p-1}(x)
\e\sum_{k=0}^{p-1}\f{\b{2k}k\b{3k}k}{(-x)^{3k}}
 \e P_{[\f p3]}\Big(1+\f{54}{x^3}\Big)
\\&\e -\Ls p3\sum_{n=0}^{p-1}\Ls {n^3-3x(x^3-216)n
-2x^6-1080x^3+108^2}p \mod p.\endalign$$
 \endpro
Proof. Since $\b{p-1}k\e (-1)^k\mod p$, using Lemma 2.1 and Fermat's
little theorem we see that
$$\sum_{k=0}^{p-1}\f{W_k(x+m)}{m^k} \e \sum_{k=0}^{p-1}\b
{p-1}kW_k(x+m)(-m)^{p-1-k}=W_{p-1}(x)\mod p.$$ On the other hand,
since $p\mid \b{2k}k\b{3k}k$ for $\f p3<k<p$ we have
$$W_{p-1}(x)=\sum_{k=0}^{[\f{p-1}3]}
\b{2k}k\b{3k}k\b{p-1}{3k}x^{p-1-3k} \e\sum_{k=0}^{p-1}
\f{\b{2k}k\b{3k}k}{(-x)^{3k}}\mod p.$$ By [28, Corollary 3.1],
$$\align \sum_{k=0}^{p-1}\f{\b{2k}k\b{3k}k}{(-x)^{3k}}
& \e P_{[\f p3]}\Big(1+\f{54}{x^3}\Big) \e -\Ls
p3\sum_{n=0}^{p-1}\Ls {n^3-3(1-\f{216}{x^3})n+
\f{108^2}{x^6}-\f{1080}{x^3}-2}p
\\&=-\Ls p3\sum_{n=0}^{p-1}\Ls {\sls n{x^2}^3
-3(1-\f{216}{x^3})\f n{x^2}+ \f{108^2}{x^6}-\f{1080}{x^3}-2}p
\\&=-\Ls p3\sum_{n=0}^{p-1}\Ls {n^3-3x(x^3-216)n
-2x^6-1080x^3+108^2}p\mod p.
\endalign$$
Thus the lemma is proved.

\par Let $p>3$ be a prime. Taking $m=1$ and $x=-4$ in Lemma 2.2 yields
$$\sum_{k=0}^{p-1}W_k\e
-\Ls{-6}p\sum_{n=0}^{p-1}\Ls{n^3-840n+9074}p\mod p.\tag 2.1$$ Taking
$m=-1$ and $x=-2$ in Lemma 2.2 yields
$$\sum_{k=0}^{p-1}(-1)^kW_k\e
-\Ls{-6}p\sum_{n=0}^{p-1}\Ls{n^3-336n+2522}p\mod p.\tag 2.2$$

\pro{Lemma 2.3} For any nonnegative integer $n$ we have
$$\sum_{k=0}^n\b nkW_k3^{n-k}=\cases \b{2n/3}{n/3}\b n{n/3}&
\t{if $3\mid n$,}
\\0&\t{if $3\nmid n$.}
\endcases$$
\endpro
Proof. Putting $x=-3$ and $y=3$ in Lemma 2.1 gives
$$\sum_{k=0}^n\b nkW_k3^{n-k}=W_n(0)=\sum_{k=0}^{[n/3]}\b{2k}k\b{3k}k\b
n{3k}0^{n-3k}.$$ This yields the result.
\par Now we are ready to prove the following result.
\pro{Theorem 2.1} Let $p$ be a prime with $p>3$. Then
$$\align&\sum_{k=0}^{p-1}\f{W_k}{(-3)^k}\e\sum_{k=0}^{p-1}\f{W_k}{(-9)^k}
\\&\e \cases -L\mod p&\t{if $p\e 1\mod 3$ and so $4p=L^2+27M^2$ with $L\e
1\mod 3$,}
\\0\mod p&\t{if $p\e 2\mod 3$}\endcases\endalign$$
and
$$\sum_{k=0}^{p-1}\b{2k}k\f{W_k}{(-12)^k}
\e \cases L^2\mod p&\t{if $p\e 1\mod 3$ and so $4p=L^2+27M^2$,}
\\0\mod p&\t{if $p\e 2\mod 3$}.\endcases$$
\endpro
Proof. Putting $m=-9$ and $x=6$ in Lemma 2.2 we get
$\sum_{k=0}^{p-1}\f{W_k}{(-9)^k}\e P_{[\f p3]}(\f 54)\mod p$. Now
applying [28, Theorem 3.2] gives the congruence for
$\sum_{k=0}^{p-1}\f{W_k}{(-9)^k}\mod p$. Since $\b{p-1}k\e
(-1)^k\mod p$ for $k=0,1,\ldots,p-1$, using Lemma 2.3 and [5,
Theorem 9.2.1] we see that
$$\align\sum_{k=0}^{p-1}\f{W_k}{(-3)^k}&\e\sum_{k=0}^{p-1}\b{p-1}kW_k3^{p-1-k}
\\&=\cases \b{\f{2(p-1)}3}{\f{p-1}3}\b{p-1}{\f{p-1}3}\e
\b{\f{2(p-1)}3}{\f{p-1}3}\e -L\mod p\\\qq\q\qq\ \t{if $3\mid p-1$
and $4p=L^2+27M^2$ with $L\e 1\mod 3$,}
\\0\mod p\q\t{if $3\mid p-2$}.\endcases\endalign$$
\par Note that $p\mid \b{2k}k$ for $k=\f{p+1}2,\ldots,p-1$
and $\b{\f{p-1}2}k\e \b{-\f 12}k=(-4)^{-k}\b{2k}k\mod p$ for
$k=0,1,\ldots,\f{p-1}2$. Using Lemma 2.3 we see that
$$\align\sum_{k=0}^{p-1}\b{2k}k\f{W_k}{(-12)^k}&\e
\sum_{k=0}^{(p-1)/2}\b{\f{p-1}2}k\f{W_k}{3^k} \e \Ls{3}p
\sum_{k=0}^{(p-1)/2}\b{\f{p-1}2}k W_k\cdot 3^{\f{p-1}2-k}
\\&=\cases \ls 3p\b{\f{p-1}3}{\f{p-1}6}\b{\f{p-1}2}{\f{p-1}6}\e 2^{\f{p-1}3}
\b{\f{p-1}2}{\f{p-1}6}^2\mod p& \t{if $3\mid p-1$,}
\\0\mod p&\t{if $3\mid p-2$}.
\endcases\endalign$$
Now assume $p\e 1\mod 3$. Then $p=A^2+3B^2$ and $4p=L^2+27M^2$ with
$A,B,L,M\in\Bbb Z$ and $A\e L\e 1\mod 3$. By [5, p.201],
$\b{\f{p-1}2}{\f{p-1}6}\e 2A\mod p$. If $2$ is a cubic residue of
$p$, then $2^{\f{p-1}3}\e 1\mod p$. It is  well known that $3\mid B$
and so $L=-2A$. Hence
$$\sum_{k=0}^{p-1}\b{2k}k\f{W_k}{(-12)^k}\e 2^{\f{p-1}3}
\b{\f{p-1}2}{\f{p-1}6}^2\e (2A)^2\e L^2\mod p.$$  Now assume that
$2$ is  a cubic nonresidue of $p$. Then $2^{\f{p-1}3}\not\e 1\mod
p$, $3\nmid B$ and $2\nmid LM$. We choose the sign of $M$ so that
$M\e L\mod 4$ and  the sign of $B$ so that $B\e A\e 1\mod 3$. From
[23, p.227] we know that
$$2^{\f{p-1}3}\e \f{-1-A/B}2\mod p,\q
A=\f{L-9M}4\qtq{and}B=\f{L+3M}4.$$ Hence
$$\align \sum_{k=0}^{p-1}\b{2k}k\f{W_k}{(-12)^k}
&\e 2^{\f{p-1}3} \b{\f{p-1}2}{\f{p-1}6}^2\e \f{-1-\f AB}2\cdot 4A^2
\e  \f{-1-\f AB}2\cdot 4(-3B^2)
\\&=6(A+B)B=6\Big(\f{L-9M}4+\f{L+3M}4\Big)\f{L+3M}4
\\&=\f 14(3L^2-27M^2)\e \f 14(3L^2+L^2)=L^2\mod p.
\endalign$$
This proves the remaining part and the proof is now complete.
\par\q
\par{\bf Remark 2.1} In [37, Conjecture 1.4] Zhi-Wei Sun
conjectured that if $p$ is a prime such that $p\e 1\mod 6$ and so
$4p=L^2+27M^2$ with $L\e 1\mod 3$, then
$$ \sum_{k=0}^{p-1}\f{W_k}{(-9)^k}\e \sum_{k=0}^{p-1}\f{W_k}{(-3)^k}
\e -L+ \f pL\mod {p^2};$$ if $p$ is a prime with $p\e 5\mod 6$, then
$\sum_{k=0}^{p-1}\f{W_k}{(-9)^k}\e 0\mod {p^2}.$
 \pro{Lemma 2.4}
Let $p$ be an odd prime, $n,x\in\Bbb Z_p$ and $n(n+4x)\not\e 0\mod
p$. Then
$$\Ls {n+4x}p\sum_{k=0}^{p-1}\b{2k}k\f{W_k(x)}{(n+4x)^k}\e \Ls{-1}p
W_{\f{p-1}2}\Big(-\f n4\Big)\e\Ls
{n}p\sum_{k=0}^{p-1}\f{\b{2k}k\b{3k}k\b{6k}{3k}} {n^{3k}}\mod p.$$
\endpro
Proof. As $\b{\f{p-1}2}k\e \b{-\f 12}k=\b{2k}k(-4)^{-k}\mod p$ and
$p\mid \b{2k}k \b{3k}k\b{6k}{3k}$ for $\f p6<k<p$, using Lemma 2.1
we see that
$$\align &\sum_{k=0}^{p-1}\b{2k}k\f{W_k(x)}{(n+4x)^k}
\\&\e \sum_{k=0}^{(p-1)/2}\b{\f{p-1}2}kW_k(x)\Ls {-4}{n+4x}^k
\e \Ls{-4(n+4x)}p\sum_{k=0}^{(p-1)/2}\b{\f{p-1}2}kW_k(x)
\Ls{n+4x}{-4}^{\f{p-1}2-k}
\\&=\Ls{-n-4x}pW_{\f{p-1}2}\Big(-\f n4\Big)
=\Ls{-n-4x}p\sum_{k=0}^{[p/6]}\b{2k}k\b{3k}k\b{\f{p-1}2}{3k}
\Big(-\f n4\Big)^{\f{p-1}2-3k}
\\&\e \Ls{n(n+4x)}p\sum_{k=0}^{[p/6]}\b{2k}k\b{3k}k\b{6k}{3k}\f
1{(-4)^{3k}\cdot (-n/4)^{3k}}
\\&\e \Ls {n(n+4x)}p\sum_{k=0}^{p-1}\f{\b{2k}k\b{3k}k\b{6k}{3k}}
{n^{3k}}\mod p.
\endalign$$
This proves the theorem.
\par\q
\pro{Theorem 2.2} Let $p$ be an odd prime. Then
$$W_{\f{p-1}2}\e \cases 4x^2\mod p&\t{if $p\e 1\mod 4$ and so
$p=x^2+y^2$ with $2\nmid x$,}
\\0\mod p&\t{if $p\e 3\mod 4$.}
\endcases$$
\endpro
Proof. Since $W_1=-3$ we see that the result is true for $p=3$. Now
assume $p>3$. Putting $n=12$ in Lemma 2.4 yields
$$W_{\f{p-1}2}=W_{\f{p-1}2}(-3)\e\Ls{-12}p\sum_{k=0}^{p-1}
\f{\b{2k}k\b{3k}k\b{6k}{3k}}{12^{3k}}\mod p.$$ In [19] Mortenson
proved the congruence
$$\sum_{k=0}^{p-1}
\f{\b{2k}k\b{3k}k\b{6k}{3k}}{12^{3k}}\e\cases \sls {-3}p4x^2\mod
p&\t{if $4\mid p-1$ and so $p=x^2+y^2$ with $2\nmid x$,}
\\0\mod p&\t{if $4\mid p-3$,}\endcases$$
which was conjectured by Rodriguez-Villegas in 2003. Now combining
the above gives the result.
\par\q
 \pro{Theorem 2.3} Let $p$ be a prime with
$p>3$. Then
$$\sum_{k=0}^{p-1}\b{2k}k\f{W_k}{54^k}
\e\cases \ls p3(4x^2-2p)\mod p&\t{if $p=x^2+4y^2\e 1\mod 4$,}
\\0\mod p&\t{if $p\e 3\mod 4$.}
\endcases$$\endpro
Proof. It is easy to check the result for $p=11$. Now assume that
$p\not=11$. Taking $n=66$ and $x=-3$ in Lemma 2.4 and then applying
[29, Theorem 4.3] deduces the result.
\par\q
\pro{Theorem 2.4} Let $p>3$ be a prime. Then
$$\sum_{k=0}^{p-1}\b{2k}k\f{W_k}{8^k}
\e\cases 4x^2-2p\mod p&\t{if $p=x^2+2y^2\e 1,3\mod 8$,}
\\0\mod p&\t{if $p\e 5,7\mod 8$.}
\endcases$$\endpro
Proof. It is easy to check the result for $p=5$. Now assume that
$p>5$. Taking $n=20$ and $x=-3$ in Lemma 2.4 and then applying [29,
Theorem 4.4] we deduce the result.
\par\q
\pro{Theorem 2.5} Let $p$ be a prime with $p\not=2,3,7$. Then
$$\align&\Ls {-3}p\sum_{k=0}^{p-1}\b{2k}k\f{W_k}{(-27)^k}
\e \Ls {-3}p\sum_{k=0}^{p-1}\b{2k}k\f{W_k}{243^k}
\\&\e\cases 4x^2-2p\mod
p&\t{if $p=x^2+7y^2\e 1,2,4\mod 7$,}
\\0\mod p&\t{if $p\e 3,5,6\mod 7$.}
\endcases\endalign$$\endpro

Proof. It is easy to check the result for $p=5,17$. Now assume that
$p\not=5,17$. Taking $n=-15$ and $x=-3$ in Lemma 2.4 and then
applying [29, Theorem 4.7] we deduce the congruence for
$\sum_{k=0}^{p-1}\b{2k}k\f{W_k}{(-27)^k}\mod p$. Taking $n=255$ and
$x=-3$ in Lemma 2.4 and then applying [29, Theorem 4.7] we deduce
the remaining part.
\par\q
\pro{Theorem 2.6} Let $p$ be a prime and $p\not=2,11$. Then
$$\align \sum_{k=0}^{p-1}\b{2k}k\f{W_k}{(-44)^k}
\e\cases x^2-2p\mod p&\t{if $\sls p{11}=1$ and so $4p=x^2+11y^2$,}
\\0\mod p&\t{if $\sls p{11}=-1$.}
\endcases\endalign$$\endpro
Proof. It is easy to check the result for $p=3$. Now assume that
$p\not=3$. Taking $n=-32$ and $x=-3$ in Lemma 2.4 and then applying
[29, Theorem 4.8] we deduce the result.
 \pro{Theorem
2.7} Let $p$ be a prime with $p\not=2,3,19$. Then
$$\align \Ls {-3}p\sum_{k=0}^{p-1}\b{2k}k\f{W_k}{(-108)^k}
\e\cases x^2-2p\mod p&\t{if $\sls p{19}=1$ and so
 $4p=x^2+19y^2$,}
\\0\mod p&\t{if $\sls p{19}=-1$.}
\endcases\endalign$$\endpro
Proof. Taking $n=-96$ and $x=-3$ in Lemma 2.4 and then applying [29,
Theorem 4.9] we deduce the result.
\par Using Lemma 2.4 and [29, Theorem 4.9] one can also deduce the
following results.

\pro{Theorem 2.8} Let $p$ be a prime, $p\not=2,3,43$. Then
$$\align \Ls {-3}p\sum_{k=0}^{p-1}\b{2k}k\f{W_k}{(-972)^k}
\e\cases x^2-2p\mod p&\t{if $\sls p{43}=1$ and so
 $4p=x^2+43y^2$,}
\\0\mod p&\t{if $\sls p{43}=-1$.}
\endcases\endalign$$\endpro

\pro{Theorem 2.9} Let $p$ be a prime with $p\not=2,3,7,67$. Then
$$\align \Ls {-3}p\sum_{k=0}^{p-1}\b{2k}k\f{W_k}{(-5292)^k}
\e\cases x^2-2p\mod p&\t{if $\sls p{67}=1$ and so
 $4p=x^2+67y^2$,}
\\0\mod p&\t{if $\sls p{67}=-1$.}
\endcases\endalign$$\endpro

\pro{Theorem 2.10} Let $p$ be a prime with  $p\not=2,3,7,11,163$.
Then
$$\align &\Ls {-3}p\sum_{k=0}^{p-1}\b{2k}k
\f{W_k}{(-640332)^k}
\\&\e\cases x^2-2p\mod p&\t{if $\sls p{163}=1$ and so
 $4p=x^2+163y^2$,}
\\0\mod p&\t{if $\sls p{163}=-1$.}
\endcases\endalign$$\endpro

\pro{Theorem 2.11} Suppose that $p>3$ is a prime, $x\in\Bbb Z_p$ and
$x(x^3+27)(x^3-216)(x^2+6x-18)\not\e 0\mod p$. Then
$$\align W_{p-1}(x)^2&\e\Big(\sum_{k=0}^{p-1}\f{W_k}{(-x-3)^k}\Big)^2
\e\sum_{k=0}^{p-1}\b{2k}k^2\b{3k}k\Big(-\f{x^3+27}{x^6}\Big)^k
\\&\e\Ls {x(x^3-216)}p
\sum_{k=0}^{p-1}\b{2k}k\b{3k}k\b{6k}{3k}\Big(-\f{x^3+27}
{x(x^3-216)}\Big)^{3k}
\\&\e
\Ls{x^3+27}pW_{\f{p-1}2}\Ls{x(x^3-216)}{4(x^3+27)}
\\&\e\sum_{k=0}^{p-1}\b{2k}k\Big(-\f{x^3+27}
{(x^2+6x-18)^2}\Big)^kW_k\mod p.\endalign$$
\endpro
Proof. From Lemma 2.2 we see that
$W_{p-1}(x)\e\sum_{k=0}^{p-1}\f{W_k}{(-x-3)^k}
\e\sum_{k=0}^{p-1}\f{\b{2k}k\b{3k}k}{(-x)^{3k}} \mod p.$ By [28,
Theorem 2.1],
$$\Big(\sum_{k=0}^{p-1}\b{2k}k\b{3k}km^k\Big)^2\e
\sum_{k=0}^{p-1}\b{2k}k^2\b{3k}k(m(1-27m))^k\mod {p^2}.$$ Hence,
$$\align& W_{p-1}(x)^2\\&\e\Big(\sum_{k=0}^{p-1}\f{W_k}{(-x-3)^k}
\Big)^2\e\Big(\sum_{k=0}^{p-1}\f{\b{2k}k\b{3k}k}{(-x)^{3k}}\Big)^2
\e \sum_{k=0}^{p-1}\b{2k}k^2\b{3k}k\Big(-\f
1{x^3}\Big(1+\f{27}{x^3}\Big)\Big)^k\mod p.\endalign$$ By [30,
Theorem 2.2], for $t\in\Bbb Z_p$ with $4t\not\e \pm 5\mod p$,
$$\sum_{k=0}^{p-1}\b{2k}k^2\b{3k}k\Ls{1-t^2}{108}^k
\e\Ls{4t+5}p\sum_{k=0}^{p-1}\b{2k}k\b{3k}k\b{6k}{3k}\Ls{(t+1)(1-t)^3}{432(4t+5)^3}^k\mod
p.$$ Taking $t=-1-\f{54}{x^3}$ gives
$$\align &\sum_{k=0}^{p-1}\b{2k}k^2\b{3k}k\Big(-\f
1{x^3}\Big(1+\f{27}{x^3}\Big)\Big)^k\\& \e \Ls {x(x^3-216)}p
\sum_{k=0}^{p-1}\b{2k}k\b{3k}k\b{6k}{3k}\Big(-\f{x^3+27}
{x(x^3-216)}\Big)^{3k}\mod p.\endalign$$ Set
$n=-\f{x(x^3-216)}{x^3+27}$. Then
$n-12=-\f{x^4+12x^3-216x+324}{x^3+27}=-\f{(x^2+6x-18)^2}{x^3+27}.$
From Lemma 2.4 (with $x=-3$) and the fact $\b{(p-1)/2}r\e
4^{-r}\b{2r}r\mod p$ for $0\le r\le \f{p-1}2$ we see that
$$\align&\sum_{k=0}^{p-1}\b{2k}k\Big(-\f{x^3+27}
{(x^2+6x-18)^2}\Big)^kW_k
\\&\e\Ls{n-12}p(-1)^{\f{p-1}2}W_{\f{p-1}2}\Big(-\f n4\Big)
=\Ls{x^3+27}pW_{\f{p-1}2}\Ls{x(x^3-216)}{4(x^3+27)}
\\&=\Ls{x^3+27}p\sum_{k=0}^{p-1}\b{2k}k\b{3k}k\b{\f{p-1}2}{3k}
\Ls{x(x^3-216)}{4(x^3+27)}^{\f{p-1}2-3k}
\\&\e\Ls {x(x^3-216)}p
\sum_{k=0}^{p-1}\b{2k}k\b{3k}k\b{6k}{3k}\Big(-\f{x^3+27}
{x(x^3-216)}\Big)^{3k}\mod p.\endalign$$ Now putting all the above
together proves the theorem.
\par\q
\pro{Corollary 2.1}  Suppose that $p>3$ is a prime, $x\in\Bbb Z_p$
and $(x+3)(1+9x+27x^2)(1+9x)(1+27x^2)(1-27x^2)\not\e 0\mod p$. Then
$$\Big(\sum_{k=0}^{p-1}W_kx^k\Big)^2
\e\sum_{k=0}^{p-1}\b{2k}k\Ls{x(1+9x+27x^2)}{(1-27x^2)^2}^kW_k\mod
p.$$
\endpro
Proof. Substituting $x$ with $-\f 1x-3$ in Theorem 2.11 yields
$$\align \Big(\sum_{k=0}^{p-1}W_kx^k\Big)^2
&\e\sum_{k=0}^{p-1}\b{2k}k\Big(-\f{-(\f 1x+3)^3+27}{(\f
1{x^2}-27)^2}\Big)^kW_k
\\&=\sum_{k=0}^{p-1}\b{2k}k\Ls{x(1+9x+27x^2)}{(1-27x^2)^2}^kW_k\mod
p.\endalign$$ This proves the corollary.
\par We remark that Corollary 2.1 is the p-analogue of (1.3).

\section*{3. Congruences for generalized Ap\'ery-like numbers}
\par In 2010, Jarvis and Verrill [13] established a relation between $u_n$
and $u_{p-1-n}$ modulo $p$ for $\{A_n'\},\{a_n\},\{f_n\}$ and
$\{S_n\}$, where $p$ is an odd prime. Inspired by (1.1) and (1.2),
we introduce generalized Ap\'ery-like numbers and prove a vast
generalization of those congruences given in [13].

 \pro{Theorem 3.1} Suppose $r\in\Bbb Z^+$ and $c\in\Bbb Z$
with $c\not=0$. Let $b(n)$ be the polynomial of $n$ with integral
coefficients and the property $b(-1-n)=(-1)^rb(n)$ for any $n\in\Bbb
Z$. Define the sequence $\{u_n\}$ by
$$u_0=1,\ u_1=b(0)\qtq{and} (n+1)^ru_{n+1}=b(n)u_n-cn^ru_{n-1}\
(n\ge 1).\tag 3.1$$ Suppose that $p$ is an odd prime with $p\nmid c$
and $u_p\in\Bbb Z_p$. For $n=0,1,2,\ldots,p-1$ we have
$$u_n\e u_{p-1}c^nu_{p-1-n}\e\cases \ls cpc^nu_{p-1-n}\mod p&\t{if $p\nmid u_{\f{p-1}2}$,}
\\ (-1)^{r-1}\ls cpc^nu_{p-1-n}\mod p&\t{if $p\mid u_{\f{p-1}2}$.}\endcases$$
In particular,
$$u_{p-1}\e \cases \ls cp\mod p&
\t{if $p\nmid u_{\f{p-1}2}$,}
\\ (-1)^{r-1}\ls cp\mod p&\t{if $p\mid u_{\f{p-1}2}$.}
\endcases$$
\endpro
Proof. By (3.1), for $n\in\{0,1,\ldots,p-1\}$, $u_n\in\Bbb Z_p$,
$(p-n)^ru_{p-n}=b(p-1-n)u_{p-1-n}-c(p-1-n)^ru_{p-2-n}$ and so
$(-n)^ru_{p-n}\e b(-1-n)u_{p-1-n}-c(-n-1)^ru_{p-2-n}\mod p.$ Since
$b(-1-n)=(-1)^rb(n)$ we get $n^ru_{p-n}\e
b(n)u_{p-1-n}-c(n+1)^ru_{p-2-n}\mod p.$
 Multiplying $c^n$ on both sides gives
$$(n+1)^rc^{n+1}u_{p-2-n}\e b(n)c^nu_{p-1-n}-cn^r\cdot
c^{n-1}u_{p-n}\mod p.\tag 3.2$$ By (3.1),
$p^ru_p=b(p-1)u_{p-1}-c(p-1)^ru_{p-2}$. Thus $b(p-1)u_{p-1}\e
c(-1)^ru_{p-2}\mod p$. Since $b(p-1)\e b(-1)=(-1)^rb(0)\mod p$ we
see that $b(0)u_{p-1}\e cu_{p-2}\mod p$.
 If $p\mid u_{p-1}$, we must have $p\mid u_{p-2}$ and
so $p\mid u_{p-3}$ by (3.1). If $u_{p-(m-1)}\e u_{p-m}\e 0\mod p$
for $m\in\{2,3,\ldots,p-1\}$, then $u_{p-(m+1)}\e 0\mod p$ by (3.1).
Hence  $u_{p-2}\e u_{p-3}\e \cdots \e u_1\e u_0\e 0\mod p$. But
$u_0=1$. This is a contradiction. Therefore $p\nmid u_{p-1}$. Set
$v_n=c^nu_{p-1-n}/u_{p-1}$. Then $v_0=1=u_0$ and
$v_1=cu_{p-2}/u_{p-1}\e b(0)=u_1\mod p$. By (3.2), for
$n=1,2,\ldots,p-1$ we have $(n+1)^rv_{n+1}\e b(n)v_n-cn^rv_{n-1}\mod
p.$ Hence $u_n\e v_n=c^nu_{p-1-n}/u_{p-1}\mod p$ for
$n=0,1,\ldots,p-1$ and so $u_{p-1}\e c^{p-1}u_0/u_{p-1}\mod p$,
which implies $u_{p-1}^2\e c^{p-1}\e 1\mod p$ and so $u_{p-1}\e
\varepsilon_p\mod p$ for $\varepsilon_p\in\{1,-1\}$. This yields
$u_n\e c^nu_{p-1-n}/u_{p-1} \e \varepsilon_pc^nu_{p-1-n}\mod p.$
Taking $n=\f{p-1}2$ gives $u_{\f{p-1}2}\e
\varepsilon_pc^{\f{p-1}2}u_{\f{p-1}2}\e \varepsilon_p\sls
cpu_{\f{p-1}2}\mod p$. Hence, if $p\nmid u_{\f{p-1}2}$, then
$\varepsilon_p\sls cp\e 1\mod p$, $\varepsilon_p=\sls cp$ and so
$u_n\e \sls cpc^nu_{p-1-n}\mod p$. Now assume that $p\mid
u_{\f{p-1}2}$. By the above argument, $u_{\f{p+1}2}\e \varepsilon_p
c^{\f{p+1}2}u_{\f{p-3}2}\e \varepsilon_p c\sls cpu_{\f{p-3}2}\mod
p$. By (3.1),
$$\Ls{p+1}2^ru_{\f{p+1}2}=b\Ls{p-1}2u_{\f{p-1}2}-c\Ls{p-1}2^ru_{\f{p-3}2}
\e -c\Ls{p-1}2^ru_{\f{p-3}2}\mod p.$$ Namely, $u_{\f{p+1}2}\e
(-1)^{r-1}cu_{\f{p-3}2}\mod p$. Hence $c(\varepsilon_p\sls
cp-(-1)^{r-1})u_{\f{p-3}2}\e u_{\f{p+1}2}-u_{\f{p+1}2}=0\mod p$. If
$p\mid u_{\f{p-3}2}$, since $p\mid u_{\f{p-1}2}$ we see that
$u_{\f{p-5}2}\e\cdots\e u_0\e 0\mod p$ by (3.1). But $u_0=1$.
Therefore $p\nmid u_{\f{p-3}2}$ and so $\varepsilon_p\sls
cp=(-1)^{r-1}$. This yields $u_n\e
\varepsilon_pc^nu_{p-1-n}=(-1)^{r-1}\sls cpc^nu_{p-1-n}\mod p$,
which completes the proof.
\par\q
 \pro{Corollary 3.1} Let $p>3$ be a
prime and $n\in\{0,1,\ldots,p-1\}$. Then
$$\align &P_n(x)\e P_{p-1-n}(x)\mod p,\q A_n\e A_{p-1-n}\mod p,
\\&D_n\e 64^nD_{p-1-n}\mod p,\q b_n\e 81^nb_{p-1-n}\mod p,\q
\\&T_n\e 16^nT_{p-1-n}\mod p,
\q W_n\e \Ls p3 27^nW_{p-1-n}\mod p,
\\&Q_n\e \Ls p3 72^nQ_{p-1-n}\mod p.\endalign$$
\endpro
Proof. Taking $u_n=P_n(x),A_n,D_n,b_n,T_n$ in Theorem 3.1 yields the
first five congruences. Since $\b{p-1}m\e (-1)^m\mod p$, using [25,
Corollary 2.2] we deduce that
$$W_{p-1}=\sum_{k=0}^{[p/3]}\b{2k}k\b{3k}k\b{p-1}{3k}(-3)^{p-1-3k}
\e\sum_{k=0}^{[p/3]}\f{\b{2k}k\b{3k}k}{27^k}\e\Ls p3\mod p.$$ Recall
that $(n+1)^2W_{n+1}=(-9n(n+1)-3)W_n-27n^2W_{n-1}\ (n\ge 1)$.
Applying Theorem 3.1 yields the result for $W_n$. Using [32, Lemma
2.4] and [25, Corollary 2.2] we see that
$$ Q_{p-1}=\sum_{k=0}^{p-1}\b{p-1}k(-8)^{p-1-k}f_k
\e\sum_{k=0}^{p-1}\f{f_k}{8^k}
\e\sum_{k=0}^{p-1}\f{\b{2k}k\b{3k}k}{27^k}\e\Ls p3\mod p.$$ Recall
that $(n+1)^2Q_{n+1}=(-17n(n+1)-6)Q_n-72n^2Q_{n-1}\ (n\ge 1)$.
Applying Theorem 3.1 yields the result for $Q_n$.
\par\q
 \pro{Theorem 3.2}
Let $\{u_n\}$ be given in Theorem 3.1, and let $p$ be an odd prime.
Suppose that $u_m\in\Bbb Z_p$ for $m=0,1,2,\ldots$ and $k\in\Bbb
Z^+$. Then
$$u_{kp+n}\e u_{kp}u_n\mod p\qtq{for}n=0,1,\ldots,p-1.$$
\endpro
Proof. We prove the theorem by induction on $n$. Clearly the result
is true for $n=0$ since $u_0=1$. By (3.1),
$(kp+1)^ru_{kp+1}=b(kp)u_{kp}-c(kp)^ru_{kp-1}$. Thus, $u_{kp+1}\e
b(kp)u_{kp}\e b(0)u_{kp}=u_1u_{kp}\mod p$. This shows that the
result is also true for $n=1$. Now assume $2\le m\le p-1$ and the
result holds for $n<m$. From (3.1) and the inductive hypothesis we
see that
$$\align m^ru_{kp+m}&\e b(kp+m-1)u_{kp+m-1}-c(kp+m-1)^ru_{kp+m-2}\\& \e
b(m-1)u_{kp+m-1}-c(m-1)^ru_{kp+m-2}
\\&\e b(m-1)u_{kp}u_{m-1}-c(m-1)^ru_{kp}u_{m-2}
=u_{kp}\cdot m^ru_m\mod p.\endalign$$ Since $p\nmid m$ we get
$u_{kp+m}\e u_{kp}u_m\mod p$. This shows that the result is true for
$n=m$. Hence the theorem is proved by induction.
\par\q
 \pro{Corollary
3.2} Let $\{u_n\}$ be given in Theorem 3.1, and let $p$ be an odd
prime. Suppose $u_m\in\Bbb Z_p$ and $u_{mp}\e u_m\mod p$ for
$m=1,2,3,\ldots$. For $n\in\Bbb Z^+$ write
$n=n_0+n_1p+\cdots+n_sp^s$, where
$n_0,n_1,\ldots,n_s\in\{0,1,\ldots,p-1\}$. Then we have the Lucas
congruence $u_n\e u_{n_0}u_{n_1}\cdots u_{n_s}\mod p$.
\endpro
Proof. Set $k=n_1+n_2p+\cdots+n_sp^{s-1}$. Then $n=kp+n_0$. By
Theorem 3.2,
$$\align u_n&\e u_{kp}u_{n_0}\e u_ku_{n_0}
=u_{n_1+n_2p+\cdots+n_sp^{s-1}}u_{n_0}\e
u_{n_0}u_{n_1}u_{(n_2+n_3p+\cdots+n_sp^{s-2})p}
\\&\e u_{n_0}u_{n_1}u_{n_2+n_3p+\cdots+n_sp^{s-2}}
\e u_{n_0}u_{n_1}u_{n_2}u_{n_3+n_4p+\cdots +n_sp^{s-3}}\\& \e\cdots
\e u_{n_0}u_{n_1}\cdots u_{n_{s-2}}u_{n_{s-1}+n_sp} \e
u_{n_0}u_{n_1}\cdots u_{n_{s-1}}u_{n_s}\mod p.\endalign$$
 This proves the corollary.
 \par\q
\par{\bf Remark 3.1} Suppose that $p$ is a prime. For $m,n\in\Bbb Z^+$ with $m\le n$
write $n=n_0+n_1p+\cdots+n_sp^s$ and $m=m_0+m_1p+\cdots+m_sp^s$,
where $n_0,\ldots,n_s,m_0,\ldots,m_s\in\{0,1,\ldots,p-1\}$. Then $\b
nm\e \b{n_0}{m_0}\b{n_1}{m_1}\cdots \b{n_s}{m_s}\mod p.$ This is
called Lucas theorem. From [13] and [17] we know that many
Ap\'ery-like numbers satisfy the Lucas congruences.
\par\q
 \pro{Theorem
3.3} Let $\{u_n\}$ be given by (3.1). Then
$$\sum_{k=0}^{n-1}b(k)(-c)^{n-1-k}u_k^2=n^ru_nu_{n-1}\q(n=1,2,3,\ldots).$$
Thus, if $p$ is an odd prime, $p\nmid c$ and $u_m\in\Bbb Z_p$ for
$m\in\Bbb Z^+$, then $$\sum_{k=0}^{p-1}\f{b(k)}{(-c)^k}u_k^2\e
0\mod{p^r}.$$
\endpro
Proof. Since
$$\f{(k+1)^ru_{k+1}}{(-c)^k}-\f{k^ru_{k-1}}{(-c)^{k-1}}
=\f{(k+1)^ru_{k+1}+ck^ru_{k-1}}{(-c)^k}=\f{b(k)u_k}{(-c)^k},$$ we
see that
$$\sum_{k=0}^{n-1}\f{b(k)}{(-c)^k}u_k^2
=\sum_{k=0}^{n-1}\Big(\f{(k+1)^ru_{k+1}u_k}{(-c)^k}-\f{k^ru_ku_{k-1}}{(-c)^{k-1}}\Big)
=\f{n^ru_nu_{n-1}}{(-c)^{n-1}}.$$ This yields the result.
\endpro
\par\q
\par As an example, taking $u_n=P_n(x)$ in Theorem 3.3 gives
$$\sum_{k=0}^{n-1}(-1)^{n-1-k}(2k+1)P_k(x)^2
=n\f{P_n(x)P_{n-1}(x)}x.\tag 3.3$$
\section*{4. Conjectures on congruences involving binomial
coefficients and Ap\'ery-like numbers}
\par The Bernoulli numbers
$\{B_n\}$, Euler numbers $\{E_n\}$ and the sequence $\{U_n\}$ are
defined by
$$\align &B_0=1,\q\sum_{k=0}^{n-1}\b nkB_k=0\q(n\ge 2),
\\& E_0=1,\q E_n=-\sum_{k=1}^{[n/2]}\b n{2k}E_{n-2k}\q(n\ge 1),
\\& U_0=1,\q U_n=-2\sum_{k=1}^{[n/2]}\b n{2k}U_{n-2k}\q(n\ge 1).
\endalign$$
For congruences involving $B_n,E_n$ and $U_n$ see [22,24,26].
\par Based on calculations with Maple, we pose the following
challenging conjectures:

\pro{Conjecture 4.1} Let $p>3$ be a prime. Then
$$\align &A_{p-1}\e 1+\f 23p^3B_{p-3}\mod {p^4},
\q D_{p-1}\e 64^{p-1}-\f{p^3}6B_{p-3}\mod {p^4},
\\&b_{p-1}\e 81^{p-1}-\f 2{27}p^3B_{p-3}\mod{p^4},
\q T_{p-1}\e 16^{p-1}+\f{p^3}4B_{p-3}\mod {p^4}.\endalign$$
\endpro
\par{\bf Remark 4.1} In [33] the author proved that
$T_{p-1}\e 16^{p-1}\mod {p^3}$ for any prime $p>3$.
\pro{Conjecture
4.2} Let $p>3$ be a prime. Then
$$\align &A'_{p-1}\e 1+\f 53p^3B_{p-3}\mod {p^4},
\  f_{p-1}\e 8^{p-1}+\f 58p^3B_{p-3}\mod {p^4},
\\& S_{p-1}\e
(-1)^{\f{p-1}2}32^{p-1}+p^2E_{p-3}\mod {p^3}, \  a_{p-1}\e\Ls
p39^{p-1}+p^2U_{p-3}\mod {p^3},
\\& W_{p-1}\e \Ls p327^{p-1}+p^2U_{p-3}\mod {p^3},
\  Q_{p-1}\e \Ls p372^{p-1}+\f 52p^2U_{p-3}\mod {p^3}.
\endalign$$
\endpro
\par{\bf Remark 4.2} Let $p>3$ be a prime. In [38] Z.W. Sun proved a congruence equivalent
to $f_{p-1}\e 8^{p-1}\mod {p^3}$. In [33] the author proved a
congruence equivalent to $S_{p-1}\e (-1)^{\f{p-1}2}32^{p-1}\mod
{p^2}$.
\par Let $p$ be an odd prime. In 2000 Ahlgren and Ono[2] proved
Beukers' conjecture $A_{\f{p-1}2}\e c(p)\mod {p^2}$, where
$\{c(n)\}$ is given by
$$\Phi_1(q)=q\prod_{k=1}^{\infty}(1-q^{2k})^4(1-q^{4k})^4=\sum_{n=1}^{\infty}
c(n)q^n\q(|q|<1).$$ It is well known that $\Phi_1(\t{\rm e}^{2\pi
iz})$ is a newform in $S_4(\Gamma_0(8))$. For $|q|<1$ define
$$\align
&\Phi_4(q)=q\prod_{k=1}^{\infty}(1-q^{4k})^6=\sum_{n=1}^{\infty}\alpha(n)q^n,
\\&\Phi_2(q)=q\prod_{k=1}^{\infty}(1-q^k)^2(1-q^{2k})(1-q^{4k})
(1-q^{8k})^2 =\sum_{n=1}^{\infty}\beta(n)q^n,
\\&\Phi_3(q)=q\prod_{k=1}^{\infty}(1-q^{2k})^3(1-q^{6k})^3=
\sum_{n=1}^{\infty}\gamma(n)q^n.
\endalign$$
It is known that $\Phi_4(\t{\rm e}^{2\pi iz})$ is a weight $3$
newform with complex multiplication by $\Bbb Q(\sqrt{-1})$, and for
$m\in\{2,3\}$ $\Phi_m(\t{\rm e}^{2\pi iz})$ is a weight $3$ newform
with complex multiplication by $\Bbb Q(\sqrt{-m})$. More precisely,
$$\align &\Phi_4(\t{\rm e}^{2\pi iz})\in
S_3\Big(\Gamma_0(16),\ls{-4}{\cdot}\Big),\\& \Phi_3(\t{\rm e}^{2\pi
iz})\in S_3\Big(\Gamma_0(12),\ls{-3}{\cdot}\Big),\q \Phi_2(\t{\rm
e}^{2\pi iz})\in S_3\Big(\Gamma_0(8),\ls{-8}{\cdot}\Big),\endalign$$
where $(\f a{\cdot})$ is the Legendre-Jacobi-Kronecker symbol. See
[11,18,21]. From [21, (14.2)] we know that for odd prime $p$,
$$\align &\alpha(p)=\cases 4x^2-2p&\t{if $p=x^2+4y^2\e 1\mod 4$},
\\0&\t{if $p\e 3\mod 4$,}\endcases\tag 4.1
\\&\beta(p)=\cases 4x^2-2p&\t{if $p=x^2+2y^2\e 1,3\mod 8$},
\\0&\t{if $p\e 5,7\mod 8$,}\endcases\tag 4.2
\\&\gamma(p)=\cases 4x^2-2p&\t{if $p=x^2+3y^2\e 1\mod 3$,}
\\0&\t{if $p\e 2\mod 3$.}\endcases\tag 4.3
\endalign$$
In [1] Ahlgren proved Beukers' conjecture:
$$A_{\f{p-1}2}'\e\cases 4x^2-2p\mod{p^2}&
\t{if $p\e 1\mod 4$ and so $p=x^2+4y^2$,}
\\0\mod{p^2}&\t{if $p>3$ and $p\e 3\mod 4$.}
\endcases$$
This is equivalent to $A'_{\f{p-1}2}\e \alpha(p)\mod {p^2}$ for
$p>3$. Using modular forms with complex multiplication,  Gomez,
McCarthy and Young [11] proved that for prime $p>2$ and $r\in\Bbb
Z^+$,
$$\align
&A'_{\f{p^r-1}2}\e\cases (x+yi)^{2r}+(x-yi)^{2r}\mod {p^2}&\t{if
$p=x^2+y^2\e 1\mod 4$ with $2\nmid x$,}
\\0\mod{p^2}&\t{if $p>3$ and $p\e 3\mod 4$,}\endcases
\\&a_{\f{p^r-1}2}\e\cases
(x+y\sqrt{-3})^{2r}+(x-y\sqrt{-3})^{2r}\mod p&\t{if $p=x^2+3y^2\e
1\mod 6$,}
\\0\mod p&\t{if $p\e 5\mod 6$,}
\endcases
\\&(-1)^{\f{p^r-1}2}f_{\f{p^r-1}2}\\&\e\cases
(x+y\sqrt{-2})^{2r}+(x-y\sqrt{-2})^{2r}\mod p&\t{if $p=x^2+2y^2\e
1,3\mod 8$,}
\\0\mod p&\t{if $p\e 5,7\mod 8$.}
\endcases\endalign$$
By [21, (13.1)], for odd prime $p$, $m\in\{1,3,5,\ldots\}$ and
$r\in\{2,3,4,\ldots\}$,
$$\aligned&A'_{\f{mp^r-1}2}
\e\cases (4x^2-2p)A'_{\f{mp^{r-1}-1}2}-p^2 A'_{\f{mp^{r-2}-1}2}
\mod{p^r}&\t{if $p=x^2+4y^2\e 1\mod 4$,}
\\p^2A'_{\f{mp^{r-2}-1}2}\mod {p^r}&\t{if $p\e 3\mod 4$,}\endcases
\\&f_{\f{mp^r-1}2}
\e\cases (-1)^{\f{p-1}2}(4x^2-2p)f_{\f{mp^{r-1}-1}2}-p^2
f_{\f{mp^{r-2}-1}2} \mod{p^r}& \t{if $p=x^2+2y^2\e 1,3\mod 8$,}
\\p^2f_{\f{mp^{r-2}-1}2}\mod {p^r}&\t{if $p\e 5,7\mod 8$,}\endcases
\\&a_{\f{mp^r-1}2}
\e\cases (4x^2-2p)a_{\f{mp^{r-1}-1}2}-p^2 a_{\f{mp^{r-2}-1}2}
\mod{p^r}&\t{if $p=x^2+3y^2\e 1\mod 3$,}
\\p^2a_{\f{mp^{r-2}-1}2}\mod {p^r}&\t{if $p\e 2\mod 3$.}\endcases
\endaligned$$
\par Now we present several conjectures, which can be viewed as generalizations
of the above results.
 \pro{Conjecture 4.3} Let $p$ be a prime of the
form $4k+3$. Then
$$\align &3\b{\f{p-3}2}{\f{p-3}4}^2A'_{\f{p-1}2}\e p^2\mod {p^3},
\\&A'_{\f{mp^r-1}2}\e p^2A'_{\f{mp^{r-2}-1}2}\mod{p^{2r}}
\q\t{for $m\in\{1,3,5,\ldots\}$ and $r\in\{2,3,4,\ldots\}$.}
\endalign$$
\endpro

\pro{Conjecture 4.4} Let $p$ be an odd prime, $m\in\{1,3,5,\ldots\}$
and $r\in\{2,3,4,\ldots\}$. Then
$$\align &a_{\f{mp^r-1}2}\e p^2a_{\f{mp^{r-2}-1}2}\mod {p^{2r-1}}
\q\t{for $p\e 5\mod 6$,}
\\&f_{\f{mp^r-1}2}\e p^2f_{\f{mp^{r-2}-1}2}\mod {p^{2r-1}}
\q\t{for $p\e 5,7\mod 8$,}
\\&S_{\f{mp^r-1}2}\e p^2S_{\f{mp^{r-2}-1}2}\mod {p^{2r-1}}
\q\t{for $p\e 5,7\mod 8$,}
\\&W_{\f{mp^r-1}2}\e p^2W_{\f{mp^{r-2}-1}2}\mod {p^{2r-1}}
\q\t{for $p\e 3\mod 4$,}
\\&Q_{\f{mp^r-1}2}\e p^2Q_{\f{mp^{r-2}-1}2}\mod {p^{2r-1}}
\q\t{for $p\e 13,17,19,23\mod {24}$.}
\endalign$$
\endpro

\pro{Conjecture 4.5} Suppose that $p$ is an odd prime,
$m\in\{1,3,5,\ldots\}$ and $r\in\{2,3,4,\ldots\}$.
\par $(\t{\rm i})$ If $p\e 1,3\mod 8$ and so $p=x^2+2y^2$, then
$$S_{\f{mp^r-1}2}\e
(4x^2-2p)S_{\f{mp^{r-1}-1}2}-p^2S_{\f{mp^{r-2}-1}2} \mod{p^r}.$$
\par $(\t{\rm ii})$ If $p\e 1\mod 4$ and so $p=x^2+4y^2$, then
$$W_{\f{mp^r-1}2}\e
(4x^2-2p)W_{\f{mp^{r-1}-1}2}-p^2W_{\f{mp^{r-2}-1}2} \mod{p^r}.$$
\par $(\t{\rm iii})$ If $p\e 1,7\mod {24}$ and so $p=x^2+6y^2$, then
$$Q_{\f{mp^r-1}2}\e
(-1)^{\f{p-1}2}(4x^2-2p)Q_{\f{mp^{r-1}-1}2}-p^2Q_{\f{mp^{r-2}-1}2}
\mod{p^r}.$$
\par $(\t{\rm iv})$ If $p\e 5,11\mod {24}$ and so $p=2x^2+3y^2$, then
$$Q_{\f{mp^r-1}2}\e
(-1)^{\f{p+1}2}(8x^2-2p)Q_{\f{mp^{r-1}-1}2}-p^2Q_{\f{mp^{r-2}-1}2}
\mod{p^r}.$$
\endpro

\pro{Conjecture 4.6} Let $\{u_n\}$ be one of the six sequences
$\{A'_n\},\ \{f_n\},\ \{S_n\},\ \{a_n\},\
 \{Q_n\}$ and $\{W_n\}$, and $c=-1,-8,32,9,72$ or $27$ according as
 $u_n=A'_n,f_n,S_n,a_n,Q_n$ or $W_n$. Suppose that $p$ is an odd prime with $p\nmid c$.
Then
$$4u_{\f{mp^2-1}2}\e (5-c^{p-1})u_{\f{p-1}2}u_{\f{mp-1}2} \mod
{p^2}\qtq{for}m=1,3,5,\ldots.$$
\endpro
\par Comparing Conjecture 4.6 with the case $r=2$ in [21,(13.1)]
 and Conjecture 4.5 suggests the following conjecture:
 \pro{Conjecture 4.7} Suppose that $p$ is an odd
prime.
\par $(\t{\rm i})$ If $p\e 1\mod 3$ and so $p=x^2+3y^2$, then
$a_{\f{p-1}2}\e (9^{p-1}+3)x^2-2p\mod{p^2}.$
\par $(\t{\rm ii})$ If $p\e 1\mod 4$ and so $p=x^2+4y^2$, then
$W_{\f{p-1}2}\e (27^{p-1}+3)x^2-2p\mod{p^2}.$
\par $(\t{\rm iii})$ If $p\e 1,7\mod {24}$ and so $p=x^2+6y^2$, then
$\sls 3pQ_{\f{p-1}2}\e (72^{p-1}+3)x^2-2p\mod{p^2}$; if $p\e
5,11\mod {24}$ and so $p=2x^2+3y^2$, then $\sls 3pQ_{\f{p-1}2}\e
(72^{p-1}+3)\cdot 2x^2-2p\mod{p^2}.$
\endpro

\par{\bf Remark 4.3} Let $p$ be an odd prime. By (4.3),
 Conjecture 4.7(i) is equivalent to
 $$a_{\f{p-1}2}\e \Big(1+\f 14(9^{p-1}-1)\Big)\gamma(p)\mod {p^2}
 \qtq{for}p\e 1\mod 3.$$
By (4.1),
 Conjecture 4.7(ii) is equivalent to
 $$W_{\f{p-1}2}\e \Big(1+\f 14(27^{p-1}-1)\Big)\alpha(p)\mod {p^2}\qtq{for}
 p\e 1\mod 4.$$
 By (4.2), the first congruence in [32, Conjecture 3.1] is equivalent to
$$f_{\f{p-1}2}\e (-1)^{\f{p-1}2}
\Big(1+\f 14(8^{p-1}-1)\Big)\beta(p)\mod {p^2}\qtq{for}
 p\e 1,3\mod 8,$$
 and the first part of [31, Conjecture 2.2] is equivalent to
$$S_{\f{p-1}2}\e
\Big(1+\f 14(32^{p-1}-1)\Big)\beta(p)\mod {p^2}\qtq{for}
 p\e 1,3\mod 8.$$

Let $p>3$ ba a prime. In [34,35] Z.W. Sun conjectured  congruences
for $\sum_{k=0}^{p-1}\b{2k}k^3/m^k$ $\mod {p^2}$ with
$m=1,-8,16,-64,256,-512,4096$. Such conjectures were proved by the
author in [27], and later proved by Kibelbek et al in [15]. In [35]
Z.W. Sun also conjectured that
$$\align&\sum_{k=0}^{p-1}\f{\b{2k}k^3}{(-8)^k}\e
\sum_{k=0}^{p-1}\f{\b{2k}k^3}{64^k}\e
(-1)^{\f{p-1}4}\sum_{k=0}^{p-1}\f{\b{2k}k^3}{(-512)^k}\mod
{p^3}\qtq{for}p\e 1\mod 4,
\\&\sum_{k=0}^{p-1}\b{2k}k^3\e
\Ls{-1}p\sum_{k=0}^{p-1}\f{\b{2k}k^3}{4096^k}\e\sum_{k=0}^{p-1}\f{\b{2k}k^2\b{4k}{2k}}{81^k}
\mod {p^3}\qtq{for}p\e 1,2,4\mod 7,
\\&\sum_{k=0}^{p-1}\f{\b{2k}k^3}{16^k}\e\Ls{-1}p\sum_{k=0}^{p-1}\f{\b{2k}k^3}{256^k}
\mod {p^3}\qtq{for}p\e 1\mod 3.
\endalign$$
In light of related calculations on Maple we made the following
conjectures.
 \pro{Conjecture 4.8} Let $p$ be a prime with $p\not=2,3,7$.
\par $(\t{\rm i})$ If $p\e 1,2,4\mod 7$ and so
$p=x^2+7y^2$, then
$$\sum_{k=0}^{p-1}\b{2k}k^3\e
(-1)^{\f{p-1}2}\sum_{k=0}^{p-1}\f{\b{2k}k^3}{4096^k} \e
4x^2-2p-\f{p^2}{4x^2}\mod {p^3}.$$
\par $(\t{\rm ii})$ If $p\e 3,5,6\mod 7$, then
$$\align &\sum_{k=0}^{p-1}\b{2k}k^3\e
\f{352}9(-1)^{\f{p-1}2}\sum_{k=0}^{p-1}\f{\b{2k}k^3}{4096^k}\\&\e\cases
-\f {11}4p^2\b{3[p/7]}{[p/7]}^{-2} \e
-11p^2\b{[3p/7]}{[p/7]}^{-2}\mod {p^3}&\t{if $p\e 3\mod 7$,}
\\-\f{99}{64}p^2\b{3[p/7]}{[p/7]}^{-2}\e
-11p^2\b{[6p/7]}{[2p/7]}^{-2}\mod {p^3}&\t{if $p\e 5\mod 7$,}
\\-\f{25}{176}p^2\b{3[p/7]}{[p/7]}^{-2}\e
-11p^2\b{[3p/7]}{[p/7]+1}^{-2}\mod {p^3}&\t{if $p\e 6\mod 7$.}
\endcases\endalign$$
\endpro
\par{\bf Remark 4.4}
From [5, Theorems 9.2.6, 12.9.8 and 12.9.9] we know that if $p$ is
an odd prime such that $p\e 1,2,4\mod 7$ and so $p=x^2+7y^2$, in
1848 Eisenstein proved that
$$\b{3[\f p7]}{[\f p7]}\e\cases 2x\mod p&\t{if $p\e 1\mod 7$ and $x\e 1\mod 7$,}
\\2x\mod p&\t{if $p\e 2\mod 7$ and $x\e 3\mod 7$,}
\\\f 25x\mod p&\t{if $p\e 4\mod 7$ and $x\e 2\mod 7$.}
\endcases$$
 \pro{Conjecture 4.9} Let $p>3$ be
a prime.
\par $(\t{\rm i})$ If $p\e 1\mod 3$ and so
$p=x^2+3y^2$, then
$$\sum_{k=0}^{p-1}\f{\b{2k}k^3}{16^k}\e
(-1)^{\f{p-1}2}\sum_{k=0}^{p-1}\f{\b{2k}k^3}{256^k} \e
4x^2-2p-\f{p^2}{4x^2}\mod {p^3}.$$

\par $(\t{\rm ii})$ If $p\e 2\mod 3$, then
$$\sum_{k=0}^{p-1}\f{\b{2k}k^3}{16^k}\e -8
(-1)^{\f{p-1}2}\sum_{k=0}^{p-1}\f{\b{2k}k^3}{256^k} \e
-p^2\b{\f{p-1}2}{\f{p-5}6}^{-2}\mod {p^3}.$$
\endpro

\pro{Conjecture 4.10} Let $p$ be an odd prime.
\par $(\t{\rm i})$ If $p\e 1\mod 4$ and so
$p=x^2+y^2$ with $2\nmid x$, then
$$\sum_{k=0}^{p-1}\f{\b{2k}k^3}{(-8)^k}\e
\sum_{k=0}^{p-1}\f{\b{2k}k^3}{64^k}\e
(-1)^{\f{p-1}4}\sum_{k=0}^{p-1}\f{\b{2k}k^3}{(-512)^k} \e
4x^2-2p-\f{p^2}{4x^2}\mod {p^3}.$$

\par $(\t{\rm ii})$ If $p\e 3\mod 4$, then
$$\sum_{k=0}^{p-1}\f{\b{2k}k^3}{(-8)^k}\e
-3\sum_{k=0}^{p-1}\f{\b{2k}k^3}{64^k}\e
6(-1)^{\f{p+1}4}\sum_{k=0}^{p-1}\f{\b{2k}k^3}{(-512)^k} \e \f
34p^2\b{\f{p-3}2}{\f{p-3}4}^{-2}\mod {p^3}.$$
\endpro

\pro{Conjecture 4.11} Let $p$ be an odd prime. Then
$$\align&\sum_{k=0}^{p-1}\f{\b{2k}k^3}{(-64)^k}
\\&\e\cases (-1)^{\f{p-1}2}(4x^2-2p-\f{p^2}{4x^2})\mod {p^3}&\t{if $p\e
1,3\mod 8$ and so $p=x^2+2y^2$,}
\\\f{p^2}3\b{[p/4]}{[p/8]}^{-2}\mod{p^3}&\t{if $p\e 5\mod 8$,}
\\\f 32p^2\b{[p/4]}{[p/8]}^{-2}\mod{p^3}&\t{if $p\e 7\mod 8$.}
\endcases\endalign$$
\endpro

\pro{Conjecture 4.12} Let $p$ be a prime with $p>3$. Then
$$\sum_{n=0}^{p-1}A_n\e \cases 4x^2-2p-\f{p^2}{4x^2}\mod {p^3}&\t{if $p\e
1,3\mod 8$ and so $p=x^2+2y^2$,}
\\ \f{17}{27}p^2\b{[p/4]}{[p/8]}^{-2}\mod{p^3}&\t{if $p\e 5\mod 8$,}
\\-\f {17}6p^2\b{[p/4]}{[p/8]}^{-2}\mod{p^3}&\t{if $p\e 7\mod 8$}
\endcases$$
and
$$\sum_{n=0}^{p-1}(-1)^nA_n\e \cases 4x^2-2p-\f{p^2}{4x^2}\mod {p^3}&\t{if $p\e
1\mod 3$ and so $p=x^2+3y^2$,}
\\ \f 54p^2\b{\f {p-1}2}{\f {p-5}6}^{-2}\mod{p^3}&\t{if $p\e 2\mod 3$.}
\endcases$$
\endpro
\par{\bf Remark 4.5} Z.W. Sun [36] conjectured the
congruences for $\sum_{n=0}^{p-1}A_n$ and
$\sum_{n=0}^{p-1}(-1)^nA_n$ modulo $p^2$ and proved the congruences
when the modulus is $p$. He also conjectured that
$$\align&
\sum_{k=0}^{p-1}A_k\e\sum_{k=0}^{p-1}\f{\b{2k}k^2\b{4k}{2k}}{256^k}
 \mod{p^3}\qtq{for}p\e 1,3\mod 8,
\\&
\sum_{k=0}^{p-1}(-1)^kA_k\e\sum_{k=0}^{p-1}\f{\b{2k}k^3}{16^k}
\mod{p^3}\qtq{for}p\e 1\mod 3.
\endalign$$
\par Let $p>3$ be a prime. In [39] Z.W. Sun conjectured the congruence for
$\sum_{n=0}^{p-1}\sum_{k=0}^n\b nk^4$ $\mod {p^2}$, and posed
conjectures on $\sum_{n=0}^{p-1}\f{D_n}{m^n}$ $\mod {p^2}$ for
$m=1,-2,4,-8,8,16,-32,64$. In [30] the author proved some
congruences for $\sum_{n=0}^{p-1}\f{D_n}{m^n}$ and
$\sum_{n=0}^{p-1}\f{b_n}{m^n}$ modulo $p$. Now we present
congruences for such sums modulo $p^3$.
 \pro{Conjecture 4.13} Let $p$ be a prime with $p\not=2,3,13,47$.
\par $(\t{\rm i})$ If $p\e
1,3\mod 8$ and so $p=x^2+2y^2$, then
$$\sum_{k=0}^{p-1}\f{\b{2k}k^2\b{3k}k}{8^k}  \e\sum_{n=0}^{p-1}b_n
\e \sum_{n=0}^{p-1}\f{b_n}{81^n}\e \sum_{n=0}^{p-1}\f{D_n}{8^n}\e
4x^2-2p-\f{p^2}{4x^2}\mod {p^3}.$$
\par $(\t{\rm ii})$ If $p\e
5,7\mod 8$, then
$$\align \sum_{k=0}^{p-1}\f{\b{2k}k^2\b{3k}k}{8^k}&\e -\f{33}{47}\sum_{n=0}^{p-1}b_n
\e -\f{99}{13}\sum_{n=0}^{p-1}\f{b_n}{81^n} \e
33\sum_{n=0}^{p-1}\f{D_n}{8^n}\\& \e\cases \f{11}9p^2\b{[
p/4]}{[p/8]}^{-2}\mod{p^3}&\t{if $p\e 5\mod 8$,}
\\-\f {11}2p^2\b{[p/4]}{[p/8]}^{-2}\mod{p^3}&\t{if $p\e 7\mod 8$.}
\endcases\endalign$$
\endpro

\pro{Conjecture 4.14} Let $p$ be a prime with $p\not=2,3,5,13$. If
$p\e 1,17,19,23\mod{30}$, then
$$\align &\sum_{k=0}^{p-1}\f{\b{2k}k^2\b{3k}k}{(-27)^k}
\e \sum_{k=0}^{p-1}\f{\b{2k}k^2\b{3k}k}{15^{3k}}\e\Ls
p3\sum_{n=0}^{p-1}\sum_{k=0}^n\b nk^4 \e \sum_{n=0}^{p-1}D_n\e
\sum_{n=0}^{p-1}\f{D_n}{64^n}
\\& \e\cases 4x^2-2p-\f{p^2}{4x^2}\mod {p^3}&\t{if $p\e
1,19\mod {30}$ and so $p=x^2+15y^2$,}
\\2p-12x^2+\f{p^2}{12x^2}\mod {p^3}&\t{if $p\e
17,23\mod {30}$ and so $p=3x^2+5y^2$.}
\endcases\endalign$$
If $p\e 7,11,13,29\mod{30}$, then
$$\align \sum_{k=0}^{p-1}\f{\b{2k}k^2\b{3k}k}{(-27)^k}
&\e 28\sum_{k=0}^{p-1}\f{\b{2k}k^2\b{3k}k}{15^{3k}} \e\Ls
p3\sum_{n=0}^{p-1}\sum_{k=0}^n\b nk^4 \e
\f{28}{53}\sum_{n=0}^{p-1}D_n\e
-\f{112}{13}\sum_{n=0}^{p-1}\f{D_n}{64^n}
\\& \e\cases
\f 72p^2\cdot 5^{[p/3]}\b{[p/3]}{[p/15]}^{-2}\mod {p^3}&\t{if $p\e
7\mod{30}$,}
\\14p^2\cdot 5^{[p/3]}\b{[p/3]}{[p/15]}^{-2}\mod {p^3}&\t{if $p\e
11\mod{30}$,}
\\ \f 7{32}p^2\cdot
5^{[p/3]}\b{[p/3]}{[p/15]}^{-2}\mod {p^3}&\t{if $p\e 13\mod{30}$,}
\\\f 78p^2\cdot 5^{[p/3]}\b{[p/3]}{[p/15]}^{-2}\mod {p^3}&\t{if $p\e
29\mod{30}$.}
\endcases\endalign$$

\pro{Conjecture 4.15} Let $p>3$ be a prime.
\par $(\t{\rm i})$ If $p\e
1\mod 3$ and so $p=x^2+3y^2$, then
$$\align \sum_{k=0}^{p-1}\f{\b{2k}k^2\b{3k}k}{108^k}
&\e\sum_{n=0}^{p-1}\b{2n}n\f{f_n}{(-4)^n}
\e\sum_{k=0}^{p-1}\f{\b{2k}k^2\b{3k}k}{1458^k} \e
\sum_{n=0}^{p-1}\f{D_n}{(-2)^n}
 \e \sum_{n=0}^{p-1}\f{D_n}{4^n}
 \e
\sum_{n=0}^{p-1}\f{D_n}{16^n}\\&\e
 \sum_{n=0}^{p-1}\f{D_n}{(-32)^n}
\e \sum_{n=0}^{p-1}\f{b_n}{(-9)^n} \e 4x^2-2p-\f{p^2}{4x^2}\mod
{p^3}.\endalign$$
\par $(\t{\rm ii})$ If $p\e
2\mod 3$, then
$$\align \sum_{k=0}^{p-1}\f{\b{2k}k^2\b{3k}k}{108^k}
&\e-\sum_{n=0}^{p-1}\b{2n}n\f{f_n}{(-4)^n} \e -\f
14\sum_{n=0}^{p-1}\f{D_n}{(-2)^n}
 \e -\sum_{n=0}^{p-1}\f{D_n}{4^n}\e
2\sum_{n=0}^{p-1}\f{D_n}{16^n}\\&\e
-\sum_{n=0}^{p-1}\f{D_n}{(-32)^n}
 \e -2\sum_{n=0}^{p-1}\f{b_n}{(-9)^n}
\e -\f {p^2}2\b{\f {p-1}2}{\f {p-5}6}^{-2}\mod{p^3}.\endalign$$
\endpro

\pro{Conjecture 4.16} Let $p>3$ be a prime.
\par $(\t{\rm i})$ If $p\e 1\mod 4$, then
$$\align &\sum_{n=0}^{p-1}\f{b_n}{(-3)^n}\e
\sum_{n=0}^{p-1}\f{b_n}{(-27)^n}\e \sum_{k=0}^{p-1}
\f{\b{2k}k^2\b{4k}{2k}}{(-12288)^k} \\& \e\cases
4x^2-2p-\f{p^2}{4x^2}\mod {p^3}&\t{if $p\e 1\mod {12}$ and so
$p=x^2+9y^2$,}
\\2p-2x^2+\f{p^2}{2x^2}\mod {p^3}&\t{if $p\e 5\mod {12}$ and so
$2p=x^2+9y^2$.}\endcases\endalign$$
\par $(\t{\rm ii})$ If $p\e 3\mod 4$, then
$$\align \sum_{n=0}^{p-1}\f{b_n}{(-3)^n}&\e -15
\sum_{n=0}^{p-1}\f{b_n}{(-27)^n}\e 10\sum_{k=0}^{p-1}
\f{\b{2k}k^2\b{4k}{2k}}{(-12288)^k}
\\& \e\cases -\f
53p^2\b{[p/3]}{[p/12]}^{-2}\mod {p^3}&\t{if $p\e 7\mod{12}$,}
\\\f 56p^2\b{[p/3]}{[p/12]}^{-2}\mod {p^3}&\t{if $p\e
11\mod{12}$.}\endcases\endalign$$
\endpro

 \pro{Conjecture 4.17} Let
$p$ be a prime with $p>7$.
\par $(\t{\rm i})$ If $p\e 1,2,4\mod 7$ and so $p=x^2+7y^2$, then
$$\align \Ls p3\sum_{n=0}^{p-1}\b{2n}n\f{W_n}{(-27)^n}
&\e \sum_{k=0}^{p-1}\f{\b{2k}k^2\b{4k}{2k}}{81^k} \e
\sum_{k=0}^{p-1}\f{\b{2k}k^2\b{4k}{2k}}{(-3969)^k}\\&\e
\Ls{-15}p\sum_{k=0}^{p-1}\f{\b{2k}k\b{3k}k\b{6k}{3k}}{(-15)^{3k}} \e
4x^2-2p-\f{p^2}{4x^2}\mod{p^3}.\endalign$$
\par $(\t{\rm ii})$ If $p\e 3,5,6\mod 7$, then
$$\align&\Ls p3\sum_{n=0}^{p-1}\b{2n}n\f{W_n}{(-27)^n} \e -\f 9{40}
\sum_{k=0}^{p-1}\f{\b{2k}k^2\b{4k}{2k}}{81^k} \e\f{45}{28}
\sum_{k=0}^{p-1}\f{\b{2k}k^2\b{4k}{2k}}{(-3969)^k}
 \\&\e -\f{375}{752}
\Ls{-15}p\sum_{k=0}^{p-1}\f{\b{2k}k\b{3k}k\b{6k}{3k}}{(-15)^{3k}}
\e\cases \f 5{16}p^2\b{3[p/7]}{[p/7]}^{-2}\mod {p^3}&\t{if $p\e
3\mod 7$,}
\\\f
{45}{256}p^2\b{3[p/7]}{[p/7]}^{-2}\mod {p^3}&\t{if $p\e 5\mod 7$,}
\\\f {125}{7744}p^2\b{3[p/7]}{[p/7]}^{-2}\mod {p^3}&\t{if $p\e 6\mod 7$}.
\endcases\endalign$$
\endpro

\pro{Conjecture 4.18} Let $p$ be a prime with $p>7$ and $p\not=71$.
\par $(\t{\rm i})$ If $p\e 1,3\mod 8$ and
so $p=x^2+2y^2$, then
$$\align \sum_{k=0}^{p-1}\f{\b{2k}k^2\b{4k}{2k}}{256^k}
&\e \sum_{k=0}^{p-1}\f{\b{2k}k^2\b{4k}{2k}}{28^{4k}} \e \Ls
{-5}p\sum_{k=0}^{p-1}\f{\b{2k}k\b{3k}k\b{6k}{3k}}{20^{3k}}\e
4x^2-2p-\f {p^2}{4x^2}\mod{p^3}.\endalign$$
\par $(\t{\rm ii})$ If $p\e 5,7\mod 8$, then
$$\align \sum_{k=0}^{p-1}\f{\b{2k}k^2\b{4k}{2k}}{256^k}
&\e -\f{441}{71}\sum_{k=0}^{p-1}\f{\b{2k}k^2\b{4k}{2k}}{28^{4k}} \e
-\f{25}7\Ls
{-5}p\sum_{k=0}^{p-1}\f{\b{2k}k\b{3k}k\b{6k}{3k}}{20^{3k}}
\\&\e\cases\f{p^2}3\b{[p/4]}{[p/8]}^{-2}\mod
{p^3}&\t{if $p\e 5\mod 8$,}
\\-\f 32p^2\b{[p/4]}{[p/8]}^{-2}\mod {p^3}&\t{if $p\e 7\mod 8$.}
\endcases\endalign$$
\endpro

\pro{Conjecture 4.19} Let $p$ be a prime with $p>5$.
\par $(\t{\rm i})$ If $p\e 1\mod 3$ and so $p=x^2+3y^2$, then
$$\sum_{k=0}^{p-1}\f{\b{2k}k^2\b{4k}{2k}}{(-144)^k}
\e \Ls p5\sum_{k=0}^{p-1}\f{\b{2k}k\b{3k}k\b{6k}{3k}}{54000^k} \e
4x^2-2p-\f{p^2}{4x^2}\mod{p^3}.$$
\par $(\t{\rm ii})$ If $p\e 2\mod 3$, then
$$\sum_{k=0}^{p-1}\f{\b{2k}k^2\b{4k}{2k}}{(-144)^k}
\e 25 \Ls p5\sum_{k=0}^{p-1}\f{\b{2k}k\b{3k}k\b{6k}{3k}}{54000^k} \e
p^2\b{(p-1)/2}{(p-5)/6}^{-2}\mod {p^3}.$$
\endpro

\pro{Conjecture 4.20} Let $p$ be a prime with $p\not=2,3,11$.
\par $(\t{\rm i})$ If $p\e 1\mod 4$ and so
$p=x^2+y^2$ with $2\nmid x$,
$$\align\Ls p3\sum_{k=0}^{p-1}\f{\b{2k}k\b{3k}k\b{6k}{3k}}{12^{3k}}
&\e \Ls {33}p\sum_{k=0}^{p-1}\f{\b{2k}k\b{3k}k\b{6k}{3k}}{66^{3k}}
\e \sum_{k=0}^{p-1}\f{\b{2k}k^2\b{4k}{2k}}{648^k}\\& \e \Ls
p3\sum_{k=0}^{p-1}\b{2k}k\f{W_k}{54^k}\e 4x^2-2p-\f{p^2}{4x^2}\mod
{p^3}.\endalign$$
\par $(\t{\rm ii})$ If $p\e 3\mod 4$, then
$$\align \Ls p3\sum_{k=0}^{p-1}\f{\b{2k}k\b{3k}k\b{6k}{3k}}{12^{3k}}
&\e\f{121}{13} \Ls
{33}p\sum_{k=0}^{p-1}\f{\b{2k}k\b{3k}k\b{6k}{3k}}{66^{3k}} \e
-3\sum_{k=0}^{p-1}\f{\b{2k}k^2\b{4k}{2k}}{648^k}
\\&\e -\f 53\Ls
p3\sum_{k=0}^{p-1}\b{2k}k\f{W_k}{54^k}\e \f
5{12}p^2\b{(p-3)/2}{(p-3)/4}^{-2}\mod {p^3}.\endalign$$
\endpro

\pro{Conjecture 4.21} Let $p$ be a prime with $p>5$.
\par $(\t{\rm i})$ If $p\e 1\mod 3$ and so $4p=L^2+27M^2$, then
$$\sum_{k=0}^{p-1}\f{\b{2k}k^2\b{3k}k}{(-192)^k}\e
\Ls{10}p\sum_{k=0}^{p-1}\f{\b{2k}k\b{3k}k\b{6k}{3k}}{(-12288000)^k}
\e L^2-2p-\f{p^2}{L^2}\mod {p^3}.$$
\par $(\t{\rm ii})$ If $p\e 2\mod 3$, then
$$\align
\sum_{k=0}^{p-1}\f{\b{2k}k^2\b{3k}k}{(-192)^k}
&\e\f{800}{161}\Ls{10}p\sum_{k=0}^{p-1}\f{\b{2k}k\b{3k}k\b{6k}{3k}}{(-12288000)^k}
\\&\e\cases  3p^2\b{[2p/3]}{[p/12]}^{-2}\e\f 34p^2\b{[2p/3]}{[p/3]}^{-2}\mod {p^3}
&\t{if $12\mid p-5,$}
\\ \f 3{49}p^2\b{[2p/3]}{[p/12]}^{-2}\e \f 34p^2\b{[2p/3]}{[p/3]}^{-2}
\mod {p^3}&\t{if $12\mid p-11$}.
\endcases\endalign$$
\endpro

\pro{Conjecture 4.22} Let $p$ be a prime with $p\not=2,11,13$.
\par $(\t{\rm i})$ If $p\e 1,3,4,5,9\mod{11}$ and so $4p=u^2+11v^2$
with $u,v\in\Bbb Z$, then
$$\sum_{k=0}^{p-1}\f{\b{2k}k^2\b{3k}k}{64^k}\e
\Ls{-2}p\sum_{k=0}^{p-1}\f{\b{2k}k\b{3k}k\b{6k}{3k}}{(-32)^{3k}} \e
u^2-2p-\f {p^2}{u^2}\mod {p^3}.$$
\par $(\t{\rm ii})$ If $p\e 2,6,7,8,10\mod{11}$  and
$f=[\f p{11}]$, then
$$\align \sum_{k=0}^{p-1}\f{\b{2k}k^2\b{3k}k}{64^k}
&\e
\f{160}{39}\Ls{-2}p\sum_{k=0}^{p-1}\f{\b{2k}k\b{3k}k\b{6k}{3k}}{(-32)^{3k}}
\\&\e\cases -p^2\Ls{5\b{4f}{2f}}{2\b{3f}f\b{6f}{3f}}^2\mod {p^3}&\t{if
$p\e 2\mod {11}$,}
\\-p^2\Ls{13\b{4f}{2f}}{30\b{3f}f\b{6f}{3f}}^2\mod
{p^3}&\t{if $p\e 6\mod {11}$,}
\\-p^2\Ls{85\b{4f}{2f}}{558\b{3f}f\b{6f}{3f}}^2\mod{p^3}&\t{if $p\e
7\mod {11}$,}
\\-p^2\Ls{7\b{4f}{2f}}{148\b{3f}f\b{6f}{3f}}^2\mod{p^3}&\t{if $p\e
8\mod {11}$,}
\\-p^2\Ls{29\b{4f}{2f}}{756\b{3f}f\b{6f}{3f}}^2\mod{p^3}&\t{if $p\e
10\mod {11}$.}
\endcases\endalign$$
\endpro

\pro{Conjecture 4.23} Let $p$ be a prime with $p\not=2,3,19$.
\par $(\t{\rm i})$ If $p\e 1,4,5,6,7,9,11,16,17\mod{19}$ and so $4p=u^2+19v^2$
with $u,v\in\Bbb Z$, then
$$\Ls{-6}p\sum_{k=0}^{p-1}\f{\b{2k}k\b{3k}k\b{6k}{3k}}{(-96)^{3k}}\e
u^2-2p-\f {p^2}{u^2}\mod {p^3}.$$
\par $(\t{\rm ii})$ If $p\e 2,3\mod{19}$  and
$f=[\f p{19}]$, then
$$ \Ls{-6}p\sum_{k=0}^{p-1}\f{\b{2k}k\b{3k}k\b{6k}{3k}}{(-96)^{3k}}
\e\cases
-\f{985}{384}p^2\Ls{\b{6f}{3f}\b{10f}{2f}}{\b{6f}f\b{10f}{3f}\b{10f}{4f}}^2
\mod{p^3}&\t{if $19\mid p-2$,}
\\-\f{197}{58080}p^2\Ls{\b{6f}{3f}\b{10f}{2f}}{\b{6f}f\b{10f}{3f}\b{10f}{4f}}^2
\mod{p^3}&\t{if $19\mid p-3$.}
\endcases$$
\endpro

\pro{Conjecture 4.24} Let $p$ be a prime with $p\not=2,3,23$.
\par $(\t{\rm i})$ If $p\e 1,5,7,11\mod{24}$, then
$$\align&\sum_{k=0}^{p-1}\f{\b{2k}k^2\b{3k}k}{216^k}
\e\Ls p3\sum_{k=0}^{p-1}\f{\b{2k}k^2\b{4k}{2k}}{48^{2k}} \e \Ls
p3\sum_{n=0}^{p-1}\f{b_n}{9^n}\e \sum_{n=0}^{p-1}\f{D_n}{(-8)^n}
\\&
\e\cases 4x^2-2p-\f{p^2}{4x^2}\mod{p^3}&\t{if $p\e 1,7\mod {24}$ and
so $p=x^2+6y^2$,}
\\8x^2-2p-\f{p^2}{8x^2}\mod {p^3}
&\t{if $p\e 5,11\mod {24}$ and so $p=2x^2+3y^2$.}
\endcases\endalign$$
\par $(\t{\rm ii})$ If $p\e 13,17,19,23\mod{24}$, then
$$\sum_{k=0}^{p-1}\f{\b{2k}k^2\b{3k}k}{216^k}
\e -7\Ls p3\sum_{k=0}^{p-1}\f{\b{2k}k^2\b{4k}{2k}}{48^{2k}} \e \f
7{23} \Ls p3\sum_{n=0}^{p-1}\f{b_n}{9^n}\e -\f 7{25}
\sum_{n=0}^{p-1}\f{D_n}{(-8)^n}\mod {p^3}.$$
\endpro
\par{\bf Remark 4.6} Suppose that $p$ is a prime such that
$p\e 1,5,7,11\mod{24}$. In [35] Z.W. Sun conjectured that
$\sum_{k=0}^{p-1}\b{2k}k^2\b{3k}k/216^k \e\ls
p3\sum_{k=0}^{p-1}\b{2k}k^2\b{4k}{2k}/48^{2k}\mod {p^3}.$

\pro{Conjecture 4.25} Let $p$ be a prime such that $p>3$ and $p\e
3,5,6\mod 7$. Then
$$\Ls p3\sum_{k=0}^{p-1}\f{\b{2k}k\b{4k}{2k}}{63^k}
\e\cases \f p{3\b{3[p/7]}{[p/7]}}\mod {p^2}&\t{if $p\e 3\mod 7$,}
\\-\f p{4\b{3[p/7]}{[p/7]}}\mod {p^2}&\t{if $p\e 5\mod 7$,}
\\\f {5p}{66\b{3[p/7]}{[p/7]}}\mod {p^2}&\t{if $p\e 6\mod 7$.}
\endcases$$
\endpro
\par{\bf Remark 4.7} For the conjecture for $\sum_{k=0}^{p-1}
\b{2k}k\b{4k}{2k}/63^k \mod {p^2}$ in the cases $p\e 1,2,4\mod 7$
see [35, Conjecture 5.14].
 In [27] the author showed that these
congruences are true when the modulus is $p$. For any prime $p>3$,
in [35] Z.W. Sun also made conjectures on
$\sum_{k=0}^{p-1}\b{2k}k\b{4k}{2k}/48^k$ and
$\sum_{k=0}^{p-1}\b{2k}k\b{4k}{2k}/72^k$ modulo $p^2$.
\par\q
\pro{Conjecture 4.26} Let $p$ be a prime with $p\e 1\mod 3$. Then
$$\sum_{k=1}^{p-1}\f{kW_k}{(-9)^k}\e 0\mod {p^2}.$$
\endpro

\pro{Conjecture 4.27} Let $p$ be a prime with $p>3$. Then
$$\align&\sum_{k=0}^{p-1}\b{2k}k\f{W_k}{(-12)^k}
\\&\e \cases L^2-2p\mod {p^2}&\t{if $p\e 1\mod 3$ and so $4p=L^2+27M^2$ with $L,M\in\Bbb Z$,}
\\0\mod {p^2}&\t{if $p\e 2\mod 3$}.\endcases\endalign$$
\endpro

\pro{Conjecture 4.28} Let $p$ be an odd prime,
 $n\in\{-640320,-5280,-960,-96,-32,$ $20,255\}$ and
$n(n-12)\not\e 0\mod p$. Then
$$\sum_{k=0}^{p-1}\b{2k}k\f{W_k}{(n-12)^k}
\e\Ls {n(n-12)}p\sum_{k=0}^{p-1}\f{\b{2k}k\b{3k}k\b{6k}{3k}}
{n^{3k}}\mod {p^2}.$$ Hence the congruences in Theorems 2.4-2.10
also hold modulo $p^2$.
\endpro

\pro{Conjecture 4.29} Let $p>3$ be a prime. Then
$$\align &\sum_{k=0}^{p-1}\b{2k}k\f{3k+1}{(-16)^k}f_k\e (-1)^{\f{p-1}2}p
+p^3E_{p-3}\mod {p^4},
\\&\sum_{k=0}^{p-1}\b{2k}k\f{7k+2}{(-27)^k}W_k\e 2\Ls {-3}p p
-4p^3U_{p-3}\mod {p^4},
\\&\sum_{k=0}^{p-1}\b{2k}k\f{7k+3}{8^k}W_k\e 3\Ls {-2}pp\mod {p^2},
\\&\sum_{k=0}^{p-1}\b{2k}k\f{7k+2}{54^k}W_k\e 2\Ls {-3}pp\mod {p^2},
\\&\sum_{k=0}^{p-1}\b{2k}k\f{14k+3}{(-44)^k}W_k\e 3\Ls {-11}pp\mod
{p^2}\qtq{for}p\not=11,
\\&\sum_{k=0}^{p-1}\b{2k}k\f{38k+7}{(-108)^k}W_k\e 7\Ls {-3}pp\mod
{p^2},
\\&\sum_{k=0}^{p-1}\b{2k}k\f{133k+26}{243^k}W_k\e 26\Ls {-3}pp\mod
{p^2},
\\&\sum_{k=0}^{p-1}\b{2k}k\f{602k+85}{(-972)^k}W_k\e 85\Ls {-3}pp\mod
{p^2},
\\&\sum_{k=0}^{p-1}\b{2k}k\f{4154k+481}{(-5292)^k}W_k\e 481\Ls {-3}pp\mod
{p^2}\q\t{for $p\not=7$.}
\endalign$$
\endpro
\par{\bf Remark 4.8} In [12] Guo proved that for any odd prime $p$,
$$\sum_{k=0}^{p-1}\b{2k}k\f{3k+1}{(-16)^k}f_k\e (-1)^{\f{p-1}2}p\mod
{p^3}.$$

$$\q$$
\leftline{Zhi-Hong Sun}\newline
\leftline{School of Mathematics and
Statistics}
\leftline{Huaiyin Normal University}
 \leftline{Huaian, Jiangsu 223300}
 \leftline{P.R. China}
\leftline{Email: zhsun@hytc.edu.cn}
 \leftline{URL: http://maths.hytc.edu.cn/szh1.htm}
$$\q$$
\end{document}